\documentclass[12pt, twoside]{article}
\usepackage{amsmath,amsthm,amssymb}
\usepackage{times}
\usepackage{enumerate}

\def\titlerunning#1{\gdef\titrun{#1}}
\makeatletter
\def\author#1{\gdef\autrun{\def\and{\unskip, }#1}\gdef\@author{#1}}
\def\address#1{{\def\and{\\\hspace*{18pt}}\renewcommand{\thefootnote}{}%
\footnote {#1}}%
\markboth{\autrun}{\titrun}}
\makeatother
\def\email#1{e-mail: #1}

\numberwithin{equation}{section}

\begin{document}

\baselineskip=17pt

\titlerunning{Sequence elimination function and the formulas of prime numbers}

\title{Sequence elimination function and the formulas of prime numbers}

\author{Ahmed Diab}

\date{}

\maketitle

\address{
\email{ahmedhamdy1401753@sci.asu.edu.eg}
}

%\subjclass{11N32, 11Y55}

\begin{abstract}
we will derive a function that eliminates any sequence of equidistant numbers from the integer numbers, then we will derive its inverse.

Then we will use the Sequence elimination function to eliminate the multiples of the prime numbers from the integer numbers to produce formulas that generate all prime numbers in some interval.

\end{abstract}

\section{Introduction}
Formulas that generate prime numbers are one of the most important prime numbers problems, which have drawn the interest of mathematicians over the centuries \cite{guy2004unsolved,weisstein2005prime,crandall2006prime,hardy1979introduction,adleman1980distinguishing,ingham1990distribution}, Much work has been done on this subject like, Euler's prime generating polynomial \cite{ribenboim2000euler,crandall2006prime}, Fermat numbers \cite{krizek200217,crandall2006prime}, Mersenne primes \cite{crandall2006prime,gillies1964three} and many other works, Each of them made his own contribution to the prime numbers, which helped to understand the prime numbers more broadly than before, but the topic remains shrouded in mystery and many unanswered questions\cite{guy2004unsolved,hardy1979introduction}.

In this paper I will focus on producing a prime generating formulas, I will explain how to derive an infinite number of these formulas, which generates all of the prime numbers in some interval see \ref{sect:int}, we will start with deriving an equation that eliminates a sequence of an equidistant numbers see \ref{sect:SEF}, then we will derive its inverse see \ref{sect:ISEF}, then we will use this equation to eliminate the multiples of prime numbers one by one from the integer numbers in a certain method see \ref{sect:formulas}, this will produce the desired formulas, at the end of this paper I will make a generalized formula to generally describe my prime numbers formulas see \ref{sect:generalization}.

\section{The sequence elimination function }

\subsection{Derivation of the sequence elimination function}\label{sect:SEF}
The sequence elimination function is a function that eliminates a sequence of equidistant numbers from the integer numbers.

These are the steps to derive the sequence elimination function:

 \begin{enumerate}
 
  \item Let $ n=\ \ 1,\ \ 2,\ \ 3,\ \ldots.\  $ ,  $  y\in\mathbb{Z} $\\
  then let
   \begin{equation}
   \frac{n}{y}=p  
   \end{equation}  
    
   $p$ is an integer number when $n$ equals $y$ or a multiple of $y$, the difference between every two successive integer values of $p$ is the same, and equals $y$, we will use these properties to obtain the elimination of a sequence of equidistant numbers.
 
\item 	Adding constant $ \frac{c}{y}\ $ to the left-hand side of eq(2.1) where $ c\in\mathbb{Z} $
\begin{center}
$\frac{n}{y}+\frac{c}{y}=p$
\end{center}
\begin{equation}
\therefore\ \frac{n+c}{y}=p
\end{equation}
This constant will shift our results $p$ corresponding to $n$ depending on the value of $c$.

\item Taking the floor to the left-hand side of eq(2.2)
\begin{equation}
\therefore\ \left\lfloor\frac{n+c}{y}\right\rfloor=p
\end{equation} 
The results will be integers, This function's mechanism is adding one to the results after $y$ numbers.

\item Adding $n$ to the left-hand side of eq(2.3)
\begin{equation}
\therefore\ n+\left\lfloor\frac{n+c}{y}\right\rfloor=p        
\end{equation}
substituting by $n$ in eq(2.4), the result is the integer numbers without a specific sequence of numbers where the number of numbers between each two successive numbers in the sequence equals $y$ and the difference between these two adjacent numbers is $y+1$ .

 Let: 
 
  $x_1$: Is the first number in the eliminated sequence.\\
       $x_2$: Is the second number in the eliminated sequence.
\begin{equation}
\therefore\ x_2-x_1=y+1
\end{equation}

 \end{enumerate}
 
 Simplifying eq(2.4)
\begin{center}
$n+\left\lfloor\frac{n+c}{y}\right\rfloor=p$
\end{center}

\begin{center}
$\left\lfloor\frac{ny+n+c}{y}\right\rfloor=p$
\end{center} 

\begin{center}
$\left\lfloor\frac{(y+1)n+c}{y}\right\rfloor=p$
\end{center}
From eq(2.5)
\begin{equation}
\left\lfloor\frac{(x_2-x_1)n+c}{(x_2-x_1)-1}\right\rfloor=p
\end{equation}\\
We will call eq(2.6) by \textbf{SEF} for (Sequence Elimination Function).

To avoid some problems while using the equation we will set the following condition  
\begin{center}
$x_{1\ }\le\ y+1$
\end{center}

Example (1): Eliminate the next sequence of numbers from the integers where the difference between the adjacent numbers is 5.
\begin{center}
\begin{tabular}{| c | c | c | c | c | c | c |}
\hline
The sequence & 3 & 8 & 13 &	18 & 23	& …\\
\hline
\end{tabular}\\
\end{center}
sol.\\
      $ x_2=8 $, $x_1=3$\\      
 Sub. in eq(2.6)
 
 \begin{equation}
 \left\lfloor\frac{5n+c}{4}\right\rfloor=p      
 \end{equation}
 
 \textbf{The method of calculating constant $\textbf{c}$:}\\
Substituting by $n$ equals the first number in the sequence the result  will be the next number $(n+1)$, and the value of the equation inside the floor in this case is an integer number so we can neglect the floor.

 Applying this to the example:\\
 Substituting by $ n=3$, $p=4$ in eq(2.7)
 \begin{center}
 $\frac{15+c}{4}=4$
 \end{center}
 
 \begin{center}
 $15+c=16$
 \end{center}
 
 \begin{center}
 $c=1$
 \end{center}
 Sub. by $c$ in eq(2.7)
 \begin{center}
 $\left\lfloor\frac{5n+1}{4}\right\rfloor=p$
 \end{center}
 Sub. by n=1,\ 2,\ 3,\ ....
 
 \begin{center}
 \begin{tabular}{| c | c | c | c | c | c | c | c | c | c | c | c |}
 \hline
 n	&1&	2&	3&	4&	5&	6&	7&	8&	9&	10&	....\\
 \hline
p	&1	&2	&4	&5	&6	&7	&9	&10	&11	&12	&....\\
\hline
 \end{tabular}
  \end{center}
  
 The sequence is eliminated from the results.
 
\subsection{Derivation of the inverse of \textbf{SEF}} \label{sect:ISEF}
\begin{enumerate}

\item Let $l=1,\ 2,\ 3,\ .... $ , $m\in\mathbb{Z}$
Then let 
\begin{equation}
\frac{l}{m}=d\ 
\end{equation}
The result is an integer numbers when $l$ equals $m$ or a multiple of it.

\item adding $\frac{g}{m}$ to the left-hand side of eq(2.8) where $g$ is an integer constant 
\begin{center}
$\frac{l}{m}+\frac{g}{m}=d$
\end{center}
\begin{equation}
\therefore\ \frac{l+g}{m}=d
\end{equation}
This will shift the results depending to the value of $g$.

\item Subtracting the left-hand side of eq(2.9) from $l$
\begin{equation}
l-\frac{l+g}{m}=d    
\end{equation}

\item Taking the ceil to the left-hand side of eq(2.10)
\begin{equation}
\left\lceil l-\frac{l+g}{m}\right\rceil=d    
\end{equation}
The results will be just integers.

There is a sequence of numbers from the inputs has the same differences between each two adjacent numbers in the sequence and this difference equals $m$, the results $\left\{d\right\}$ of this sequence is the same as the results of the next input after each number in the sequence like, the result of one equals one and the result of two equals one also, this is the same for all the numbers in the sequence, this will stand to infinity as we see from the derivation.The results of the inputs without this sequence is the numbers $1, 2, 3, 4, …$, so if we took this inputs without this specific sequence we will get  the inverse effect that \textbf{SEF} makes. we will use this by adjusting our inputs and the other constants to get the inverse of \textbf{SEF}, first thing is to take the inputs $ l=p $ where $p$ is the output of \textbf{SEF} and it is the numbers $(1, 2, 3, 4, …)$ without a specific sequence of numbers with the same differences between each two adjacent numbers, this difference equals $ y+1 $, from this $m$ must equal $ y+1 $ and this is the second substitution, the last step is substituting by $ c=g $ because we want to get the same shift at the results that $c$ makes, by these adjustments the results $d$ will equal $n=1,\ 2,\ 3,\ ....$ .\\
Before substituting we will simplify eq(2.11). 

\begin{center}
$\left\lceil l-\frac{l+g}{m}\right\rceil=d$
\end{center}

\begin{center}
$\left\lceil\frac{ml-l\ -\ g}{m}\right\rceil=d$
\end{center}

\begin{equation}
\left\lceil\frac{(m-1)l-g}{m}\right\rceil=d 
\end{equation}

\item Applying the next substitutions to eq(2.12)

\begin{itemize}

\item $l=p$

\item $m=y+1$

\item $c=g$

\item $ d=n=1,\ 2,\ 3,\ ....$

\end{itemize}

So eq (2.12) will be

\begin{center}
$\left\lceil\frac{(y+1-1)p-c}{y+1}\right\rceil=d$
\end{center}

\begin{center}
$\left\lceil\frac{yp-c}{y+1}\right\rceil=n$
\end{center}

\begin{equation}
\left\lceil\frac{(x_2-x_1-1)p-c}{(x_2-x_1)}\right\rceil=n
\end{equation}

Eq(2.13) is the inverse of \textbf{SEF} and we can call it \textbf{ISEF}.

\end{enumerate}
 Example (2): solve:
 \begin{center}
 $\left\lfloor\frac{5n+1}{4}\right\rfloor=p$
 \end{center}
 sol.
 \begin{center}
$ \left\lceil\frac{4p-1}{5}\right\rceil=n$
 \end{center}

\section{Derivation of prime numbers formulas} \label{sect:formulas}

We will use eq(2.6) to make the prime numbers formulas by eliminating the multiples of prime numbers one by one, The elimination process includes the prime number itself.

\subsection{Elimination of the multiples of prime number 2}

The first prime number whose multiples we want to eliminate is 2. 

This is the produced formula after eliminating the multiples of 2 by using eq(2.6) 
\begin{equation}
2n-1=p     
\end{equation}  
where $n=\ 1,\ \ 2,\ \ 3,\ \ldots.\  $

The first number in the results is one, To shift our results and get rid of number one which is not a prime number, we will replace $n$ by $n+1$ in eq(3.1).
\begin{center}
$2(n+1)-1=p$
\end{center}
\begin{equation}
2n+1=p       
\end{equation}
where $n=\ 1,\ \ 2,\ \ 3,\ \ldots.\  $

This formula - eq(3.2) - gives us all the prime numbers in the interval $\left]2,3^2\right[$. 
 
 \subsection{Elimination of the multiples of prime number 3}
The second prime whose multiples we want to eliminate is 3, In the results of eq(3.2) we don't have all the multiples of 3 this because there are common multiples between 2 and 3 that have already been eliminated, In this case this will not create any problem (I will explain the possible problems and how to solve them while eliminating the multiples of 5) because multiples of three in the results of eq(3.2) are a sequence of equidistant numbers, So we can easily eliminate them by \textbf{SEF}.

We will eliminate the multiples of 3 in $p$ from the results of eq(3.2) by eliminating its corresponding $n$ which is $(n=1, 4, 7, 10, ....)$ from $n$, 

Using \textbf{SEF} to eliminate the sequence $(n=1, 4, 7, 10, ....)$ from $n$, 
\begin{equation}
\therefore\left\lfloor\frac{3n_1+1}{2}\right\rfloor=n
\end{equation}
Sub. by $n$ from eq(3.3) in eq(3.2)
\begin{equation}
2\left\lfloor\frac{3n_1+1}{2}\right\rfloor+1=p
\end{equation}
where $n_1=1,\ 2,\ 3,\ ....$

This is the formula that eliminates the multiples of 2 and 3 from the integers, And gives us all the prime numbers in the interval $\left]3,5^2\right[$.

\subsection{Elimination of the multiples of prime number 5}
The third prime whose multiples we want to eliminate from eq(3.4) is 5, we can easily determine it from our last formula.
 
In the results of eq(3.4) We don't have all the multiples of 5 because of the common multiples between 2 and 3 with 5 that are already eliminated, In this case the multiples remaining are not a sequence of equidistant numbers, So we cannot eliminate them directly.
 
There are common multiples between 2, 3 and 5 that you can get by multiplying the multiples of 6 by 5, Since the multiples of 6 are common multiples between 2 and 3.

The common multiples between 2, 3 and 5 are a sequence of equidistant numbers, The distribution of the multiples of 5,3 and 2 is the same after each number in this sequence, and the uncommon multiples of 5 with 2 and 3 too, From this we will divide the multiples of 5 into a number of sequences equal to the number of uncommon multiples of 5 with 2 and 3 in a region bounded by two common multiples between 2, 3 and 5, We can use the region between 0 and the first common multiple between 2, 3 and 5, The resulting sequences are a sequences of equidistant numbers, All sequences have the same difference between its adjacent numbers, So that we can eliminate its corresponding inputs from our last formula - eq(3.4) in our example -  by \textbf{SEF}, for 5 there will be two sequences because the number of the uncommon multiples of 5 with 2 and 3 from 0 to the first common multiple between 2, 3 and 5 is 2.

\subsubsection{Eliminating the first sequence}
The first number in the first sequence is 5, Its corresponding $n_1$ is 1, The second number in the first sequence is 35, its corresponding $n_1$ is 11, The corresponding number can be determined using the inverse of eq(3.4).

using \textbf{SEF} 
\begin{equation}
\left\lfloor\frac{10n_2+8}{9}\right\rfloor=n_1 
\end{equation}  
Sub. by $n_1$ in eq(3.4)
\begin{equation}
2\left\lfloor\frac{3\left\lfloor\frac{10n_2+8}{9}\right\rfloor+1}{2}\right\rfloor+1=p  
\end{equation}
Where $n_2=1,\ 2,\ 3,\ ....$

\subsubsection{Eliminating the second sequence}
The second sequence had the same differences between its adjacent numbers as the first, but the difference is reduced by one,  because we will eliminate it from the last formula - eq(3.6) -, where the first sequence is eliminated this causes a missing number between each two adjacent numbers in the second sequence,so the difference between their corresponding $n_2$ is 9.

The first number in the second sequence is 25, its corresponding $n_2$ is 7.

Using \textbf{SEF}
\begin{equation}
\left\lfloor\frac{9n_3+1}{8}\right\rfloor=n_2
\end{equation}
Sub. by $n_2$ in eq(3.6)
\begin{equation}
2\left\lfloor\frac{3\left\lfloor\frac{10\left\lfloor\frac{9n_3+1}{8}\right\rfloor+8}{9}\right\rfloor+1}{2}\right\rfloor+1=p  
\end{equation}
Where $n_3=1,\ 2,\ 3,\ ....$

Eq(3.8) gives us all the prime numbers in the interval $\left]5,7^2\right[$.

\subsection{Elimination of the multiples of prime number 7}

The fourth prime is 7, using the same steps we used with 5 we can eliminate the multiples of 7 and all the next primes, The difference will be the number of sequences into which multiples (not eliminated multiples) of this prime will be divided, for 7 it will be 8 sequences, and the formula that eliminates the multiples of 7 from eq(3.8) is:

\begin{center}
$\left\lfloor{\frac{\mathbf{50}\ast\left\lfloor{\frac{\mathbf{49}\ast\mathbf{n_4}+\mathbf{1}}{\mathbf{48}}}\right\rfloor+\mathbf{13}}{\mathbf{49}}}\right\rfloor=\mathbf{k}$
\end{center}

\begin{center}
$\left\lfloor{\frac{\mathbf{54}\ast\left\lfloor{\frac{\mathbf{53}\ast\left\lfloor{\frac{\mathbf{52}\ast\left\lfloor{\frac{\mathbf{51}\ast\mathbf{k}+\mathbf{20}}{\mathbf{50}}}\right\rfloor+\mathbf{24}}{\mathbf{51}}}\right\rfloor+\mathbf{31}}{\mathbf{52}}}\right\rfloor+\mathbf{35}}{\mathbf{53}}}\right\rfloor=\mathbf{b}$
\end{center}
\small
{\begin{equation}
\mathbf{2}\ast\left\lfloor{\frac{\mathbf{3}\ast\left\lfloor{\frac{\mathbf{10}\ast\left\lfloor{\frac{\mathbf{9}\ast\left\lfloor{\frac{\mathbf{56}\ast\left\lfloor{\frac{\mathbf{55}\ast\mathbf{b}+\mathbf{42}}{\mathbf{54}}}\right\rfloor+\mathbf{54}}{\mathbf{55}}}\right\rfloor+\mathbf{1}}{\mathbf{8}}}\right\rfloor+\mathbf{8}}{\mathbf{9}}}\right\rfloor+\mathbf{1}}{\mathbf{2}}}\right\rfloor+\mathbf{1}=\mathbf{p}
\end{equation}
\normalsize Where $n_4=1,\ 2,\ 3,\ ....$\\

This formula - eq(3.9) - gives us all the prime numbers in the interval $\left]7,11^2\right[$, I divided it into 3 parts because the length of the floor brackets will be too long for one part.

\subsection{The interval of prime numbers produced} \label{sect:int}

All the previous prime numbers formulas give us all the prime numbers in an open interval.
\begin{enumerate}

 \item The first term in this open interval is the last prime number whose multiples we have eliminated in the formula, because it has already been eliminated and also all the prime numbers smaller than it.
 \item The last term in the interval is the square of the first prime larger than the last prime eliminated, as it is the first uncommon multiple of this prime with the primes eliminated.
  \end{enumerate}

\subsection{Generalization of my prime numbers formulas} \label{sect:generalization}
This is a generalized formula that describes generally my prime numbers formulas, to be able to refer to any of these formulas easily.

\begin{equation}
{A(P_j)}^M=P
\end{equation}
The inverse is 
\begin{equation}
{A^{-1}(P_j)}^M=P
\end{equation} 
Where\\
$P_j$: The last prime number whose multiples we have eliminated.\\
$M$: The number of sequences eliminated.

These are generalized formulas for some of the prime numbers formulas in this paper: 
\begin{enumerate}

\item Eq(3.2):
\begin{center}
${A(2)}^1=P$
\end{center}

\item Eq(3.4):
\begin{center}
${A(3)}^2=P$
\end{center}

\item Eq(3.6):
\begin{center}
${A(5)}^3=P$
\end{center}

\item Eq(3.8):
\begin{center}
${A(5)}^4=P$
\end{center}

\end{enumerate}

\section{conclusion}

In this paper we have derived the sequence elimination function (\textbf{SEF}) step by step as well as its inverse, and illustrated how to use it in making the formulas of prime numbers, in a process of eliminating the multiples of prime numbers one by one from the integer numbers, we started by eliminating the multiples of 2 then the multiples of 3 then the multiples of 5, in the process of eliminating the multiples of 5 we had some problems then we showed how to solve these problems, after solving these problems we will not have any problem with eliminating the multiples of any prime greater than 5, then we showed the formula that eliminates the multiples of 7 and all prime numbers smaller than 7 which generates all prime numbers in the interval $\left]7,11^2\right[$, and we showed how to determine these intervals for all my prime numbers formulas, and finally I made a generalized formula for my prime numbers formulas.

\bibliographystyle{plain}
\bibliography{references}}

\end{document}